%% file: Building.tex
\documentclass[letterpaper, conference]{IEEEtran}  
\pdfoutput=1                                                        

\IEEEoverridecommandlockouts                              

\usepackage{graphics} 
\usepackage{epsfig} 
\usepackage{amsmath} 
\usepackage{amssymb}  
\usepackage{color}
\usepackage{float}
\usepackage{algorithm}
\usepackage{algorithmic}
\usepackage{hyperref}
\usepackage{url}
\usepackage{calrsfs}
\usepackage{textcomp}
\DeclareMathAlphabet{\pazocal}{OMS}{zplm}{m}{n}

\DeclareMathOperator*{\argmin}{arg\,min}  

\begin{document}

\title{Modeling and Optimization of Complex Building Energy Systems with Deep Neural Networks}

\author{
	\IEEEauthorblockN{Yize Chen, Yuanyuan Shi and Baosen Zhang\IEEEauthorrefmark{1}}
	\IEEEauthorblockA{Department of Electrical Engineering, University of Washington, Seattle, WA, USA
		\\\{yizechen, yyshi, zhangbao\}@uw.edu}
}

\maketitle

\begin{abstract}
Modern buildings encompass complex dynamics of multiple electrical, mechanical, and control systems. One of the biggest hurdles in applying conventional model-based optimization and control methods to building energy management is the huge cost and effort of capturing diverse and temporally correlated dynamics. Here we propose an alternative approach which is model-free and data-driven. By utilizing high volume of data coming from advanced sensors, we train a deep Recurrent Neural Networks~(RNN) which could accurately represent the operation's temporal dynamics of building complexes. The trained network is then directly fitted into a constrained optimization problem with finite horizons. By reformulating the constrained optimization as an unconstrained optimization problem, we use iterative gradient descents method with momentum to find optimal control inputs. Simulation results demonstrate proposed method's improved performances over model-based approach on both building system modeling and control.
\end{abstract}

\begin{IEEEkeywords}
	Building energy management, deep learning, gradient algorithms, HVAC systems
\end{IEEEkeywords}

\section{Introduction}
\input{intro}

\section{Problem Formulation \& Preliminaries}
\label{preliminary}
\input{preliminary}

\section{Recurrent Neural Networks}
\label{RNN_paragraph}
\input{RNN}

\section{Inputs Optimization for Building Control}
\label{method}
\input{methods}

\section{Case Study}
\label{results}
\input{results}

\section{Conclusion}
In this work, we are exploiting Recurrent Neural Networks' ability of learning complex temporal interactions among high-dimensional building dynamics. Our proposed method consists Recurrent Neural Networks regression and sequence optimization steps, which could both be solved efficiently. Our proposed approach is easily to be deployed for any building unit provided with rich historical running data. Simulation results show that our method outperforms existing ones both in capturing the thermal dynamics of the building as well as providing effective control solutions.

\bibliographystyle{IEEEtran}
\bibliography{bib}

\end{document}

%% file: intro.tex

According to a recent United Nations Environment Programme~(UNEP) report, buildings are responsible for 40\% of the global energy consumption~\cite{cheng2008kyoto}. Consequently, managing the energy consumption of buildings has significant economical, social, and environmental impacts, and has received much attention from researchers. Many approaches have been proposed to control building systems~(e.g., commercial and office buildings, data centers) for energy efficiency, such as nonlinear adaptive control, Model Predictive Control (MPC) and decentralized control for building heating, ventilation, and air conditioning (HVAC) systems~\cite{ma2012predictive, zhang2017decentralized, shi2016leveraging}. However, most previous research on building energy management are either based on the detailed physics model of buildings~\cite{trvcka2010overview} or simplified RC circuit models~\cite{ ma2012predictive, zhang2017decentralized,ma2009model}. The former often involves tedious and complex modeling processes with a huge number of variables and parameters, whereas the latter cannot fully capture the long term dynamics of large commercial buildings.

With the advance of sensing, communication and computing, detailed operation data are being collected for many buildings.
These data along with future weather forecasts can be utilized for data-driven real-time optimization approaches. In~\cite{jain2017data}, the authors developed a data predictive control method to replace the traditional MPC controller by using data to build a regression tree that represent the dynamical model for a building. However, regression trees still results in a linear model that can be far away from the true dynamics of building systems.  While in~\cite{yang2015reinforcement,Wei2017}, reinforcement learning was proposed to learn control policies without any explicit modeling assumptions, but computational costs for searching through large state and action spaces is hight. Ill-defined reward functions~(e.g., sparse, noisy and delayed rewards) could also prevent reinforcement learning algorithm finding the optimal control solutions~\cite{mnih2013playing}. Furthermore, large commercial buildings may have quality of service constraints that prevent the deep exploration of some states in reinforcement learning.

\begin{figure}[h!]
	\centering
	\includegraphics[scale=0.52]{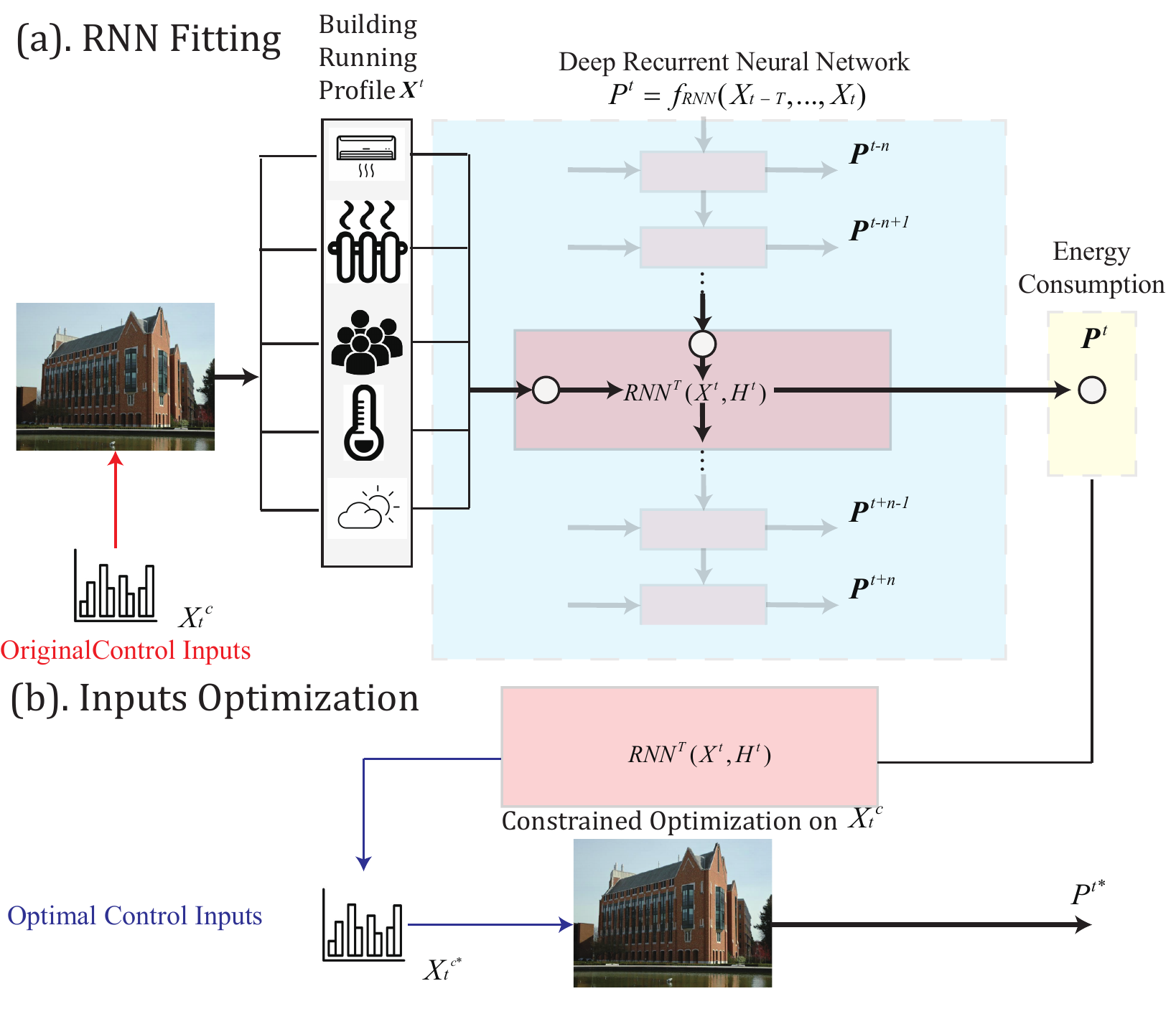}
	\caption{\small
		{
			Our model architecture for building energy system modeling and optimization based on a deep Recurrent Neural Networks~(RNN).
		}
	}
	\label{intro}
\end{figure}

In this work, we address these challenges by proposing a data-driven method which closes the loop for accurate predictive model and real-time control. The method is based on \emph{deep recurrent neural networks} that leverage rich volumes of sensor data \cite{funahashi1993approximation}. Though neural network has previously been adopted as an approach for designing controllers, the lack of large datasets and computation capabilities have prevented it from being deployed in real-time applications~\cite{nguyen1990neural}.
Firstly in a supervised learning manner, our Recurrent Neural Networks~(RNN) firstly learns the complex temporal dynamics mapping from various measurements of building operation profiles to energy consumption. Next we formulate an optimization problem with the objective of minimizing building energy consumption, which is subject to RNN-modeled building dynamics as well as physical constraints over a finite horizon of time. To solve the constrained optimization problem in a block-splitting approach, we take iterative gradient descent steps on the set of controllable inputs (e.g., zone temperature setpoints, heat rejected/added into each zone) at the current timestep. It thus finds the control inputs for each timestep. Fig. \ref{intro} illustrates our model framework. Our approach does not need analysis on complex interactions within conduction, convection or radiation processes. In addition, it can be easily scaled up to large buildings and distributed algorithms.

The main contributions of our paper are:
\begin{itemize}
	\item We model the building energy dynamics using recurrent neural networks, which leverages large volumes of data to represent the complex dynamics of buildings.
	\item We propose an input/output optimization algorithm which efficiently find the optimal control inputs for the model represented by RNN.
	\item The proposed modeling and optimization approaches open door to the integration of complex system dynamics modeling and decision-making.
\end{itemize}

The contents of the paper are as follows. The rolling horizon control problem formulation and model-based method are firstly presented in Section.~\ref{preliminary}. In Section.~\ref{RNN_paragraph} we show the design of a deep RNN which models the dynamics of complex building systems. We then reformulate the control problem as an unconstrained optimization problem, and propose the algorithm to find optimal control inputs in Section.~\ref{method}. Finally, simulation results on large building HVAC system are evaluated and compared with model-based control method in Section.~\ref{results}.

%% file: preliminary.tex
\subsection{Problem Formulation}
We consider a building energy system which includes several subsystems and zones with potentially complex interactions between them. No information about the exact system dynamics is known. At time $t$, we are provided with the building's running profile $\mathbf X_t := [\mathbf X_t^{uc}, \mathbf X_t^c, \mathbf X_t^{phy}]^T$, where $\mathbf X_t^{uc}$ denotes a collection of uncontrollable measurements such as zone temperature measurements, system node temperature measurements, lighting schedule, in-room appliances schedule, room occupancies and etc; and $\mathbf X_t^c$ denotes a collection of controllable measurements such as zone temperature setpoints, appliances working schedule and etc; $\mathbf X_t^{phy}$ denotes the set of physical measurements or forecasts values, such as dry bulb temperature, humidity and radiation volume. There are some physical constraints on some of $\mathbf X_t^c$ and $\mathbf X_t^{uc}$, for example the temperature setpoints as well as real measurements should not fall out of users' comfort regions. Without loss of generality, we denote the constraints as $ \underline{\mathbf X}_t^c \leq \mathbf X_t^c\leq \overline{\mathbf X}_t^c$ and $ \underline{\mathbf X}_t^{uc} \leq \mathbf X_t^{uc} \leq \overline{\mathbf X}_t^{uc}$. Building System operators have a group of past running profile $\mathbf{X}=\{\mathbf{X}_t\}$ along with the collection of energy consumption metering at each time step $\mathbf{P}=\{P_t\}$.

We are interested in firstly learning a model $f(\mathbf X_{t-T},...,\mathbf{X}_t)=P_t$, where $f(\cdot)$ denotes the predictive model with known parameters representing building's physical dynamics. $f(\cdot)$ maps past $T$ timestep's running profile to energy consumption at timestep $t$.

With a model $f(\cdot)$ representing the building dynamics, we formulate an optimal finite-horizon predictive control problem, and propose an efficient algorithm to find the group of optimal control inputs $\mathbf X_t^{c*}$. At timestep $t$, the control input $\mathbf{X}_{t}^c$ minimizes the energy consumption of the building for future $T$ steps. Meanwhile, previous $T$ steps' control inputs would affect current energy consumption. The objective of the controller is to minimize the energy consumption with a rolling horizon $T$, while maintaining some variables within comfortable intervals. Mathematically, we formulate the general control problem as
\begin{subequations}
	\label{formulation}
 \begin{align}
\underset{\mathbf X_t^c, ...,\mathbf X_{t+T}^c}{\text{minimize}} \quad &    
 \sum_{\tau=0}^{T}P_{t+\tau}^2 \label{1a}\\
\text{subject to} \quad 
 & P_{t+\tau}=f(\mathbf X_{t-T+\tau},...,\mathbf{X}_{t+\tau}), \forall \tau \label{1b}\\
 & \underline{\mathbf X}_{t+\tau}^c \leq \mathbf X_{t+\tau}^c\leq \overline{\mathbf X}_{t+\tau}^c, \forall \tau \label{1c}\\
 & \underline{\mathbf X}_{t+\tau}^{uc} \leq \mathbf X_{t+\tau}^{uc}\leq \overline{\mathbf X}_{t+\tau}^{uc}, \forall \tau \label{1d}\\
  \begin{split}
   & \mathbf X_{t+\tau}^{uc}=h(\mathbf X_{t-T+\tau},...,\mathbf X_{t-1+\tau},\mathbf X_{t+\tau}^c, \mathbf X_{t+\tau}^{phy}), \\ 
   &\forall \tau \label{1e}
   \end{split}
 \end{align}
\end{subequations}
where \eqref{1b} $h(\cdot)$ denotes the rolling horizon predictive model; \eqref{1c} and \eqref{1d} are the constraints on controllable and uncontrollable variables respectively; the $h(\cdot)$ in \eqref{1e} denotes a rolling horizon predictive function for uncontrollable variables based on past $T$ steps' observations as well as current step control inputs and physical forecasts. 

\subsection{First-Order Thermal Dynamic Model}
For building HVAC system, one popular method used in finite-horizon MPC to model the thermal dynamics is the reduced Resistance-Capacitance (RC) model \cite{ma2012predictive, zhang2017decentralized, ma2009model}. Here we use a rolling horizon MPC controller as a benchmark for comparison.

\begin{figure}[h!]
	\centering
	\includegraphics[scale=0.4]{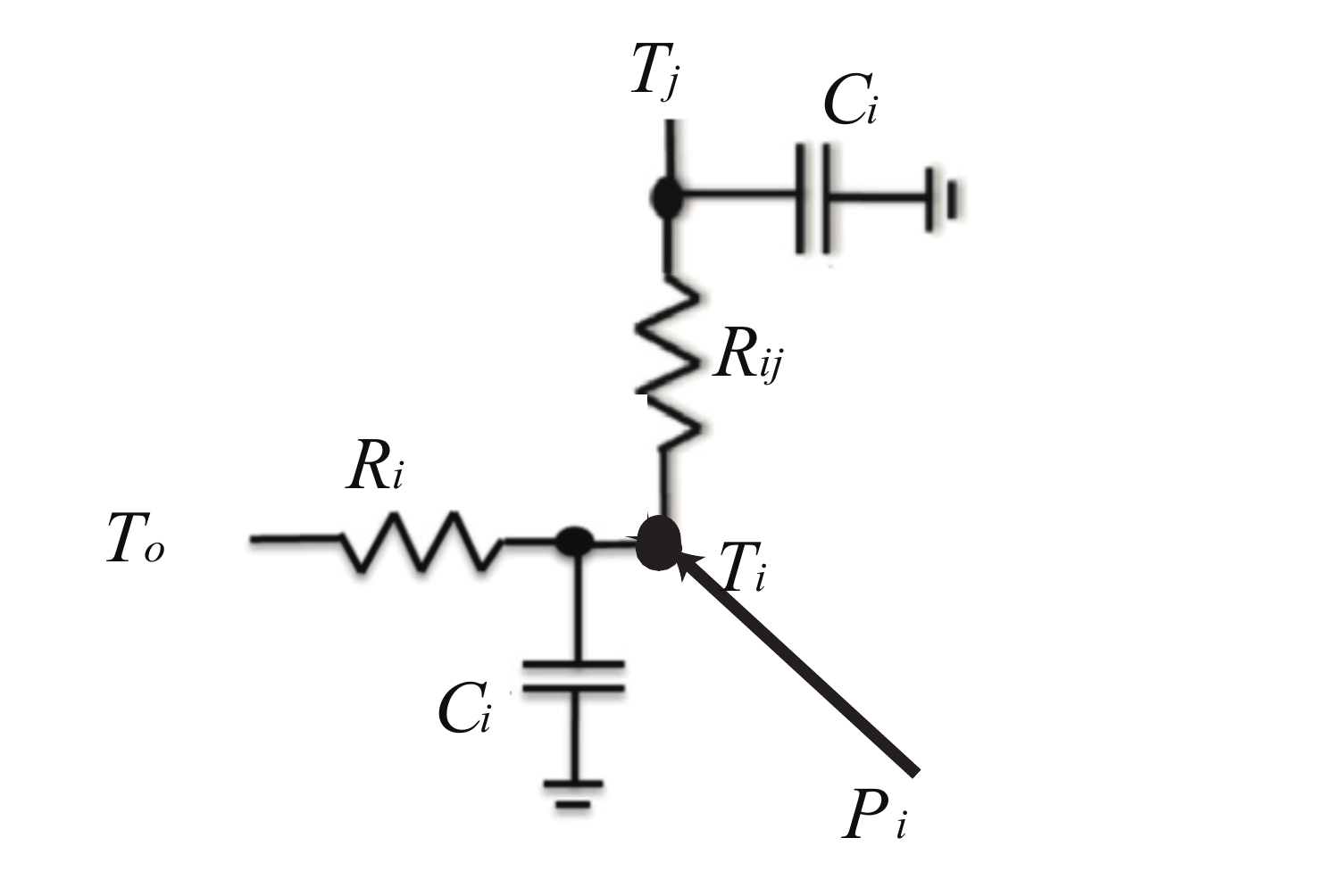}
	\caption{\small
		{
			RC network with thermal exchange between different comp.
		}
	}
	\label{RC_circuit}
\end{figure}

Denote $\pazocal{N}(i)$ as the neighboring zones for zone $i$, the first-order RC model modeling HVAC dynamics is formulated as
\begin{equation}
\label{RC}
C_i\dot{T}_{i,t}=\frac{T_{o,t}-T_{i,t}}{R_i}+\sum_{j\in \pazocal{N}(i)}^{}\frac{T_{j,t}-T_{i,t}}{R_{ij}}+P_{i,t}
\end{equation}
where $C_i, T_i$ are the thermal capacitance and room temperature for each zone $i$, while $T_o$ is the outside dry bulb temperature, and $R_i, R_{ij}$ are the thermal resistance for zone $i$ against the outside and the neighboring zone $j$. The schematic of RC network for modeling HVAC system is shown in Fig.~\ref{RC_circuit}.

Once we find $C_i, R_i, R_{ij}$ for all the zones, we have a 1st-order system to model the thermal dynamics. Since $T_i \in \mathbf X^{uc}, T_o \in \mathbf X^{phy}$, by reformulating \eqref{RC} and taking a sum of $P_{i}$ for all zones, we reformulate and write the building overall thermal dynamics
\begin{equation}
\label{RC_equ}
P_{t}=f_{RC}(\mathbf X_{t-T},...,\mathbf{X}_{t})
\end{equation}
which is further used in the optimal control problem defined in \eqref{1a}-\eqref{1e}.
MPC for building HVAC system under different model settings has been implemented in \cite{ma2012predictive,zhang2017decentralized}. We focus on the performance comparison of RC model to our proposed method in both model fitting and optimization tasks.

%% file: RNN.tex
Since the 1st-order thermal dynamic model defined in \eqref{RC} does not either capture complicated nonlinear dynamics, nor model the long-term temporal dependencies of building HVAC system, the deep RNN model becomes a good replacement.

\begin{figure}[h!]
	\centering
	\includegraphics[scale=0.43]{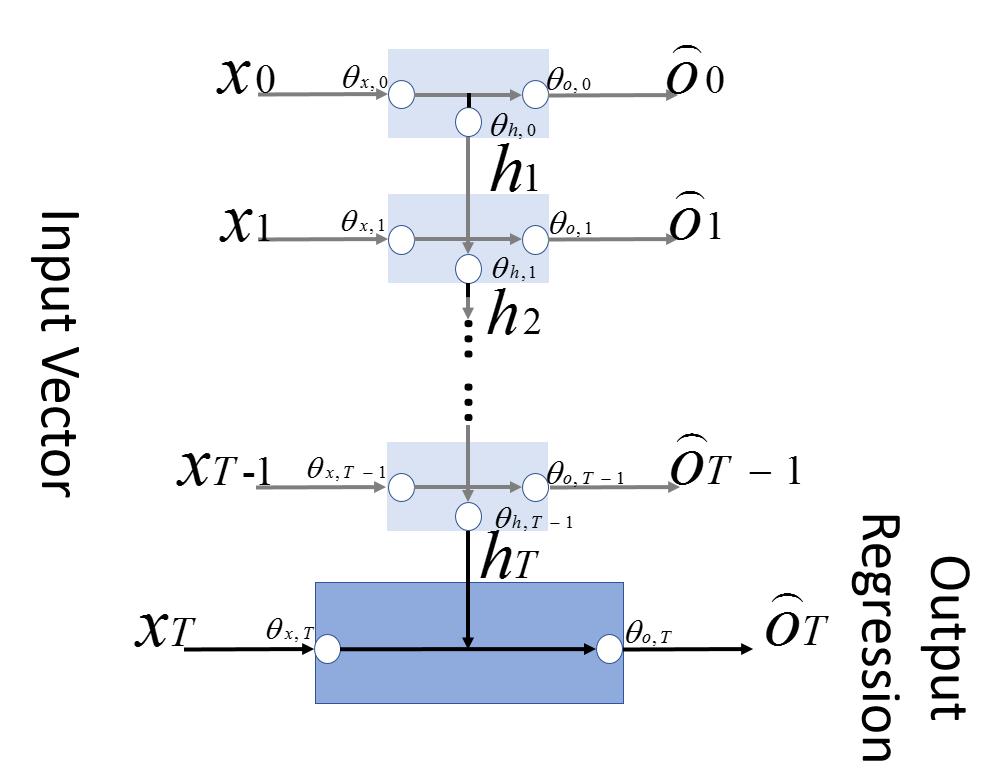}
	\caption{\small
		{
			A graphical model illustrating the RNN which is used for modeling $T$-length input-output sequential data. $\theta_{h,t},\theta_{o,t},\theta_{x,t}$ are the neural weights associated with hidden states $h_{t+1}$, output $\hat{o}_t$, and input $ x_t$ respectively.
		}
	}
	\label{RNN}
\end{figure}

RNN is a class of artificial neural networks specially designed for sequential data modeling. Unlike fully-connected neural networks where inputs are fed into the neural networks as a full vector, RNN feeds input sequentially into a neural network with directed connections. It uses its internal memory to process time-series inputs. In Fig.~\ref{RNN} we show the structure of an RNN model. 

We specifically design the RNN model to solve a time-dependent regression problem. That is to say, we want RNN automatically learn the relationship between sequential input $x_t,\: t=0,,...,T$ and output $o_T$. At timestep $t$, RNN is provided with hidden state vector $h_t$ and input vector $x_t$, and outputs its computation vector $\hat{o}_t$. The $t$-step RNN cell is composed of three group of neurons, $\theta_{x,t}, \theta_{h,t}, \theta_{o,t}$. They are associated with input, hidden state and output respectively, and are organized in function $f_{\theta_{x,t}, \: \theta_{o,t}}, f_{\theta_{h,t}, \theta_{x,t}}$ to complete the following computations:
\begin{subequations}
\begin{align}
\hat{o}_t=f_{\theta_{x,t}, \theta_{o,t}}(x_t, h_t),\\
h_{t+1}=f_{\theta_{h,t},\theta_{x,t}}(x_t, h_t)
\end{align}
\end{subequations}
where $\hat{o}_t$ is the RNN's prediction output, while $h_{t+1}$ is passed into next neuron group and takes part in $t+1$ step's computation. 

After concatenating all the neurons cells from $0$ to $T$, we get the chain function to compute $h_t$. Thus the RNN compute the final prediction value $\hat o_T$:
\begin{equation}
	\hat{o}_T=f_{\theta_{x,T}, \: \theta_{o,T}}(x_T, h_T)
\end{equation}

Since $h_t$ captures information from past inputs $x_{t-1}$, we trace hidden states back into functions of previous steps' hidden states and inputs. Thus final output $\hat{o}_T$ is eventually a function of the sequential inputs $x_t, \: t=0,...,T$. For simplicity, let's denote $\theta=\{\theta_{h,t},\theta_{o,t}, \theta_{x,t}\},\: t=0,...,T$ to be the set of neurons used in modeling the $T$-length temporal data, and wrap up all neural-composed functions of $\{f_{\theta_{x,t}, \theta_{o,t}}, f_{\theta_{h,t},\theta_{x,t}}\}$ to get the overall function $f_{RNN}$, which utilizes $\theta$ to find the output predictions with length $T$ time-series input:
\begin{equation}
\hat{o}_T=f_{RNN}(x_0,..., x_T)
\end{equation}

We set up the RNN model and initialize neuron weights $\theta$ by sampling from a normal distribution. During batch-training process, with a group of sequential input $x_t, t=0,...,T$, $\hat{o}_T$ is firstly computed, and by doing back-propagation using stochastic gradient descent~(SGD) with respect to all neurons~\cite{bottou2010large}, $\theta$ is optimized to minimize the regression loss defined in mean-square-error~(MSE) form:
\begin{subequations}
		\label{RNN_para}
	\begin{align}
	\label{RNN_loss}
	&L_{training}(\theta)=||\hat{o}_T - o_T||_2^2,\\
	\label{RNN_theta}
	&\theta^*=\argmin L_{training}(\theta)
	\end{align}

\end{subequations}

We then set up a length-$T$ RNN accordingly for our building dynamics modeling problem. With the training sets of input vectors of  historical building operating profiles $\{\mathbf X_{t-T},...,\mathbf{X}_{t}\}$ and an output energy consumption $P_t$, our RNN model $f_{RNN}$ is trained to represent the system dynamics
\begin{equation}
\label{RNN_equ}
\hat{P}_{t}=f_{RNN}(\mathbf X_{t-T},...,\mathbf{X}_{t})
\end{equation}

Our RNN is totally data-driven, and can process and represent temporal dependencies. With a rich volume of historical building operating data $\mathbf X$ and $\mathbf P$ provided as the training datasets, we train a deep RNN, which accurately models the nonlinear, complex temporal dynamics of building system. We will show in Section~\ref{results} that our deep RNN model outperforms RC model in fitting the dynamics of a large-scale building HVAC system.

%% file: methods.tex

In this section we describe our control algorithm which is based on our pre-trained deep learning model. We demonstrate how it is able to incorporate \eqref{RNN_equ} into the optimization problem \eqref{formulation}. We also illustrate how to solve such optimization problem to find a collection of optimal control sequential inputs.

By substituting $f(\cdot)$ in \eqref{formulation} with $f_{RNN}$, and denote $\mathbf {X}_t^{var}=[\mathbf {X}_t^c, \mathbf {X}_t^{uc}]$,the finite horizon control problem for building energy management is written as
\begin{subequations}
\label{rnn_formulation}
\begin{align}
\underset{\mathbf X_t^c, ...,\mathbf X_{t+T}^c}{\text{minimize}} \quad &  \sum_{\tau=0}^{T}P_{t+\tau}^2 \\
\text{subject to} \quad 
& P_{t+\tau}=f_{RNN}(\mathbf X_{t-T+\tau},...,\mathbf{X}_{t+\tau}), \forall \tau\\
& \underline{\mathbf X}_{t+\tau}^{var} \leq \mathbf X_{t+\tau}^{var}\leq \overline{\mathbf X}_{t+\tau}^{var},\forall \tau \label{9c}\\
\begin{split}
&\mathbf X_{t+\tau}^{uc}=h(\mathbf X_{t-T+\tau},...,\mathbf X_{t-1+\tau},\mathbf X_{t+\tau}^c, \mathbf X_{t+\tau}^{phy}), \\ 
&\forall \tau \label{9d}
\end{split}
\end{align}
\end{subequations}

Since $X_{t+\tau}^{uc}, \: \tau=1,...,T$ is directly controlled by control inputs of previous time. For all the uncontrollable variables with constraints we model, they also possess pairing controllable variables, e.g., the temperature measurements-temperature setpoints. We then choose $X_{t+\tau}^{uc}=X_{t-1+\tau}^{c}, \: \tau=1,...,T$, since such uncontrollable values are the control outputs corresponding to the previous step's control inputs. Thus we diminish constraint \eqref{9d}.

Since the constrained optimization problem \eqref{rnn_formulation} includes a non-convex deep neural network in the constraints, we use log barriers functions to rewrite the problem in an unconstrained form: 
\begin{equation}
\label{log_barrier}
\begin{split}
\underset{\mathbf X_t^c, ...,\mathbf X_{t+T}^c}{\text{min}} \; &L_{opti}(\mathbf X_t^c, ...,\mathbf X_{t+T}^c)=\\ &\sum_{\tau=0}^{T}f_{RNN}^2(\mathbf X_{t-T+\tau},...,\mathbf{X}_{t+\tau})\\
& -\lambda \sum_{\tau=0}^{T}\log(\mathbf X_{t+\tau}^{var}-\underline{\mathbf X}_{t+\tau}^{var})\\
&-\lambda\sum_{\tau=0}^{T}\log(\overline{\mathbf X}_{t+\tau}^{var}-\mathbf X_{t+\tau}^{var})
\end{split}
\end{equation}
where $\lambda$ is a tuning parameter, and $L_{opti}(\mathbf X_t^c, ...,\mathbf X_{t+T}^c)$ defines a loss function with inputs $\mathbf X_t^c, ...,\mathbf X_{t+T}^c$. We solve this loss minimization problem by iteratively taking gradient descents of (\ref{log_barrier}). Note that during RNN model training, we are taking gradients $\nabla_{\theta}L_{training}(\theta)$ with respect to all the neurons. Once training is done, $L_{training}(\theta)$ is converged. The RNN model serves as the temporal physical model, and is always modeling the building system dynamics accurately. Here we are taking gradients with this fixed, pre-trained RNN model, and find gradients $\nabla_{\mathbf X_{t+\tau}^c}L_{opti}(\mathbf X_t^c, ...,\mathbf X_{t+T}^c), \tau=0,...,T$ with respect to the group of controllable variables. Once $L_{opti}(\mathbf X_t^c, ...,\mathbf X_{t+T}^c)$ is converged, and we find $\mathbf {X}_t^{c^*}$ that is a local optimal solution. $\mathbf {X}_t^{c^*}$ is also the solution of controllable inputs for the finite horizon optimal control problem at timestep $t$.

The $k$-step gradient descent method is working as follows:

\begin{subequations}
\label{gradient_descent}
	\begin{align}
&g_{t+\tau,k}=\eta \nabla_{{X}_{t+\tau,k}^c}L_{opti}(\mathbf X_{t, k-1}^c, ...,\mathbf X_{t+T,k-1}^c)\\
&\mathbf{X}_{t+\tau,k}^c=\mathbf{X}_{t+\tau, k-1}^c-g_{t+\tau,k}, \tau=0,...,T
\end{align}
\end{subequations}
where $\eta$ is the learning rate, and ${X}_{t+\tau,k}^c$ denotes the value for ${X}_{t+\tau}^c$ after $k$ step's update.

Throughout our modeling and optimization approach, we do not make any physical model assumptions, and directly utilize a deep RNN to extract the model dynamics as well as finding the optimal actions to take at each time step to cut down energy consumption. We summarize the proposed method in Algorithm \ref{algorithm}, which closes the loop for building dynamics modeling and control inputs optimization. In our implementation, we improve the algorithm performance by adding momentum to gradient descents~(MomentumGD), which is shown to get over some local minima during optimization iterations as well as accelerating the convergence~\cite{sutskever2013importance}. The MomentumGD is realized as follows:

\begin{subequations}
	\label{momentumGD}
	\begin{align}
	&g_{t+\tau,k}=\gamma g_{t,k-1}+\eta \nabla_{{X}_{t+\tau,k}^c}L_{opti}(\mathbf X_{t, k-1}^c, ...,\mathbf X_{t+T,k-1}^c)\\
	&\mathbf{X}_{t+\tau,k}^c=\mathbf{X}_{t+\tau, k-1}^c-g_{t+\tau,k}, \tau=0,...,T
	\end{align}
\end{subequations}
where $\gamma$ is a momentum term determining how much previous gradients are incorporated into current step's update.

 \begin{algorithm}
	\caption{Input Optimization for Building Control}
	\label{algorithm}
	\begin{algorithmic}
		\REQUIRE Pre-trained RNN $f_{RNN}$, learning rate $\eta$, momentum $\gamma$, input optimization iterations $N_{iter}$
		\REQUIRE Control window-size $T$
		\REQUIRE Sensor measurements $\mathbf {X}_t^{uc}$, weather forecasts $\mathbf {X}_t^{phy}$
		\ENSURE $\mathbf X_t, ...,\mathbf X_{t+T}$
		\ENSURE Optimal control inputs  $X_t^*\leftarrow \emptyset$
		\FOR {iteration$=0,...,N_{iter}$}
		\STATE Update $\mathbf{X}_t^c$ using gradient descent:
		\FOR {$\tau = 0,...,T$}
		\STATE $g_{t+\tau}\leftarrow \nabla_{\mathbf{X}_{t+\tau}^c}L_{opti}(\mathbf X_t^c, ...,\mathbf X_{t+T}^c)$
		\STATE $\mathbf{X}_{t+\tau}^c \leftarrow \mathbf{X}_{t+\tau}^c - \eta \cdot MomentumGD(\mathbf{X}_{t+\tau}^c, g_{t+\tau},\gamma)$
		\ENDFOR
		\STATE Update $\mathbf{X}_t^{uc}$ using gradient descent:
		\FOR {$\tau = 0,...,T$}
		\STATE $\mathbf{X}_{t+\tau}^{uc}=\mathbf{X}_{t+\tau-1}^{c}$
		\ENDFOR
		\ENDFOR
		\STATE $\mathbf X_t^*$.insert($\mathbf{X}_{t}^c$)
	\end{algorithmic}
\end{algorithm}

%% file: results.tex
In this section, we set up a realistic model in standard building simulation software EnergyPlus \cite{crawley2001energyplus}. We demonstrate the effectiveness of our data-driven approach for both system dynamics modeling and building energy management. In order to compare with the model-based approach, we focus on the HVAC system for a large building complex. But our method is a general regression and optimization approach, which could be easily applied to overall building energy management problem.

\subsection{Experimental Setup}
\begin{figure}[h!]
	\centering
	\includegraphics[scale=0.36]{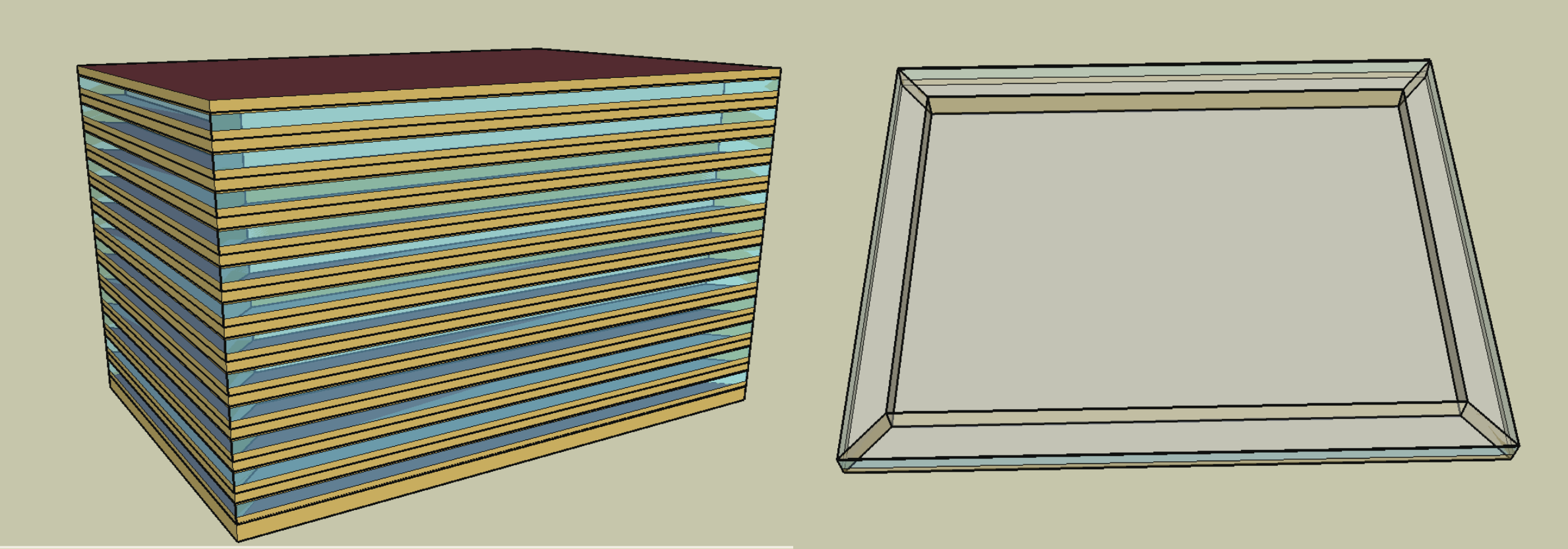}
	\caption{\small
		{
			Schematic diagram of simulated large commercial building.
		}
	}
	\label{building}
\end{figure}

We set up our EnergyPlus simulations using a 12-storey large office building (in Fig.~\ref{building}) listed in the commercial reference buildings from U.S. Department of Energy (DoE CRB)~\cite{deru2011us}. The building has a total floor area of $498,584$ square feet which is divided into $16$ separate zones. We simulate the building running through the year of $2004$ in Seattle, WA, and record ($\mathbf{X}_t$, $P_t$) with a resolution of $10$ minutes. We shuffle and separate $2$ months' data as our stand-alone testing dataset for both regression and control performance evaluation, while the remaining $10$ months' data is used to for RNN training. The processed datasets have $55$ input features, which include controllable variables such as zone temperature setpoints, and uncontrollable variables such as zone occupancies and temperature measurements. Output is a single feature for energy consumption at each timestep. We directly use historical weather data records into both RC model and RNN model. For future work, the forecasts model should also be considered into the pipeline. A finite horizon of $4$ hours is set for both MPC method and proposed method.

We set up our deep learning model using Tensorflow, a Python open-source package. Our RNN model is composed of $1$ recurrent layer with $3$ subsequent fully-connected layers. We adopt rectified linear unit~(ReLU) activation functions, dropout layers and Stochastic Gradient Descent~(SGD) optimizer to improve our neural network training.

\subsection{Simulation Results}
\begin{figure}[h!]
	\centering
	\includegraphics[scale=0.6]{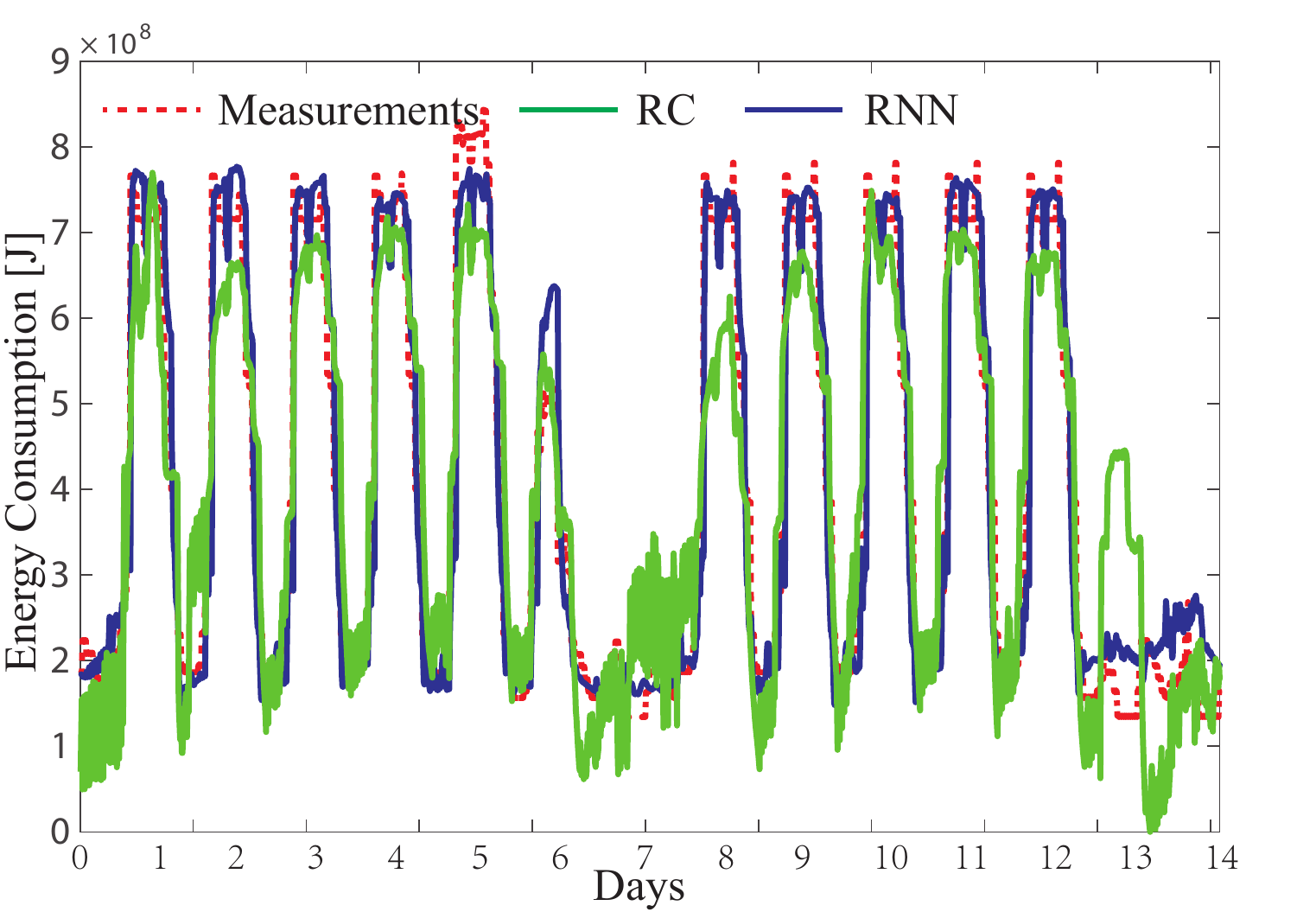}
	\caption{\small
		{
		Comparison of building's real energy consumption measurements(real), Recurrent Neural Networks' predicted energy consumption (blue) and RC model prediction (green) on a week of testing data.
		}
	}
	\label{fitting}
\end{figure}

We first compare the model fitting performance for 1st order model and RNN model, and the fitting result for two weeks' energy consumption is shown in Fig.~\ref{fitting}. To quantitatively compare the model fitting error, we calculate the Root-Mean-Square-Error~(RMSE) value for normalized energy consumption on test dataset. RMSE for the first-order RC model is $0.240$. The RNN model improves RC model by $68.33\%$ with an RMSE of $0.076$. It is also interesting to notice that this large office building actually has an energy consumption dropdown on weekdays' noon due to the occupancy schedule. RC model fails to capture this dynamic characteristics, while RNN  model is able to fit noon values given past $4$ hours' input measurements. Moreover, RC model performs poorly on weekend regression task, which hardly represents the HVAC dynamics. This inaccurate model would make subsequent MPC algorithm fail to operate on correct model space.

\begin{figure}[h!]
	\centering
	\includegraphics[scale=0.6]{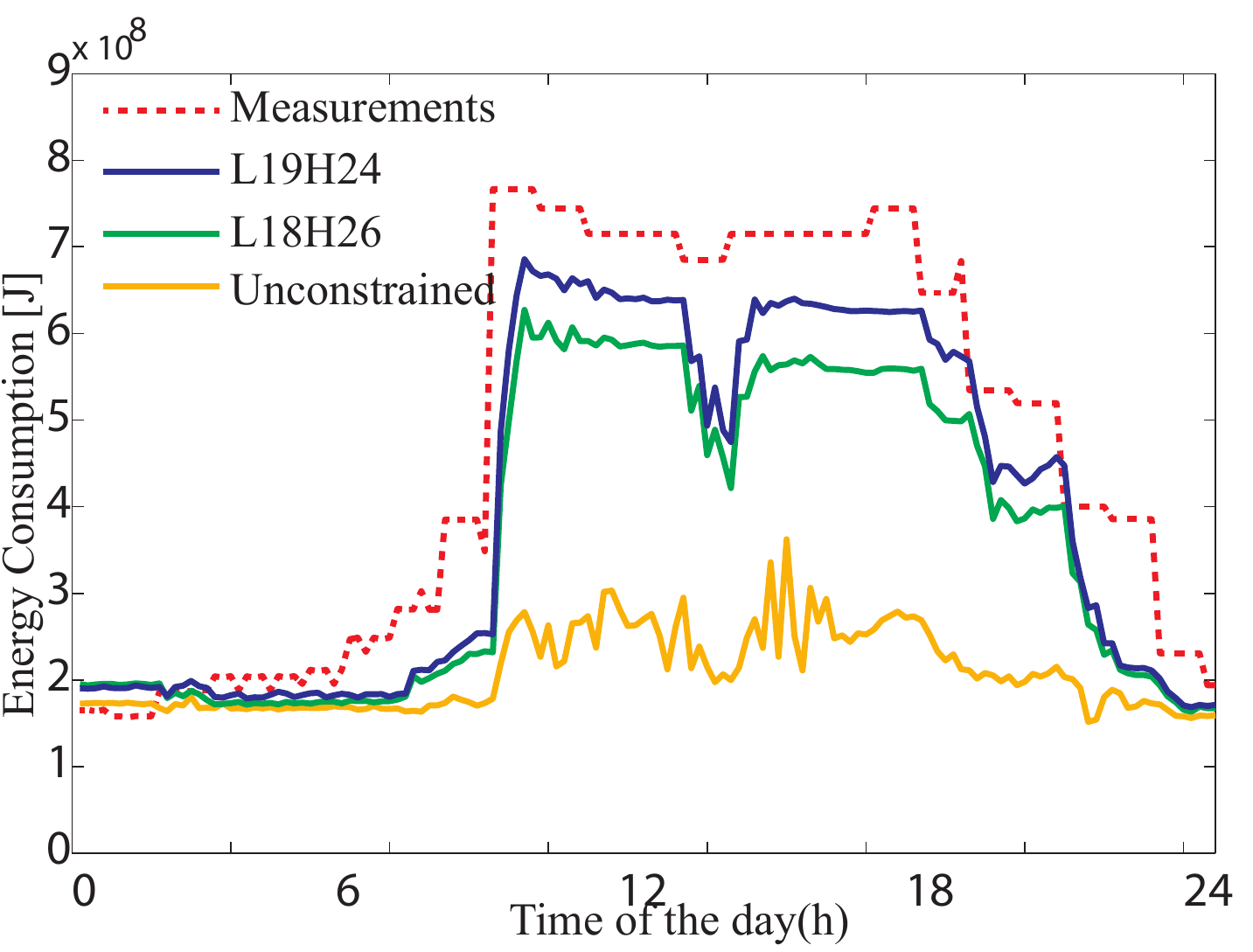}
	\caption{\small
		{
			Effects of constraints interval on optimization performance.
		}
	}
	\label{constraints}
\end{figure}

Next we show the constrained optimization problem formulated in \eqref{log_barrier} is efficient in finding optimal inputs $\mathbf{X}_t^c$ for the HVAC system. In Fig.~\ref{constraints} we show a group of 3 plots corresponding to different zone temperature setpoint constraints. We keep setpoint constraints the same for all the 16 zones. Compare the results of $\mathbf{X}_t^c \in $[18\textdegree{}C,26\textdegree{}C] and the results of $\mathbf{X}_t^c \in $[19\textdegree{}C,24\textdegree{}C], we observe that our approach is able to find sharper control inputs with less energy consumption when constraint intervals are bigger. When there is no constraint on temperature, our approach simply finds extreme control inputs such that the energy consumption is nearly same as the midnight consumption. 

\begin{figure}[h!]
\centering
\includegraphics[scale=0.55]{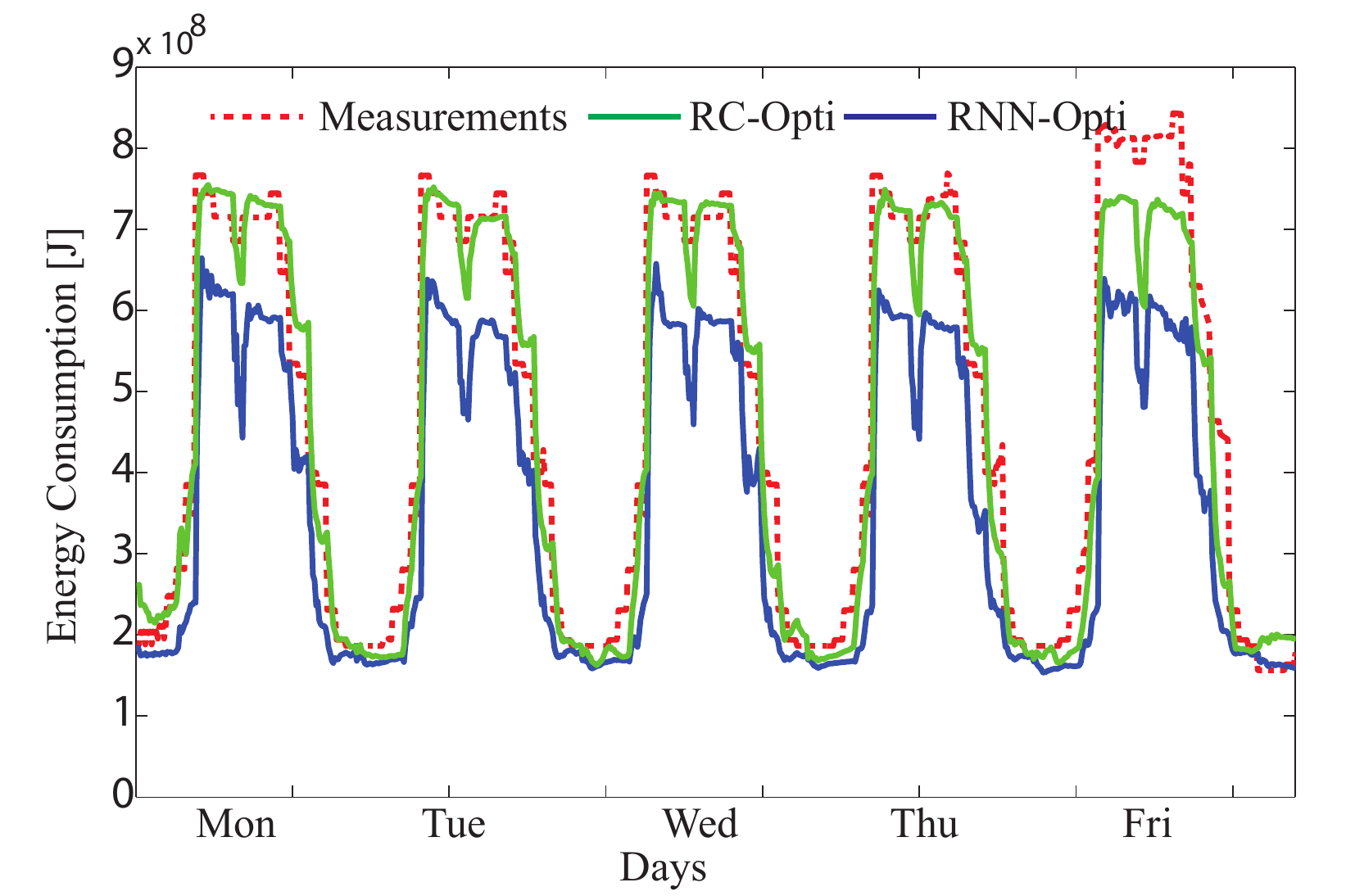}
\caption{\small
{
Comparison of optimization results of RC model~(green) and RNN model~(blue) with respect to original measurements~(red).
  }
}
\label{opti_results}
\end{figure}

We then compare the optimization performance for RC model and RNN model. Fig.~\ref{opti_results} illustrates a Monday-Friday energy consumption profile with temperature setpoint constraints $\mathbf{X}_t^c \in $[18\textdegree{}C,26\textdegree{}C]. By using RNN model and taking the gradient steps, we find a sequence of control inputs that could reduce $30.74\%$ of energy consumption. On the other hand, the solution found by RC model only gives us a $4.07\%$ reduction of energy consumption. This furthur illustrates that RC model is not good at modeling large-scale building system dynamics.

\begin{figure}[h!]
	\centering
	\includegraphics[scale=0.67]{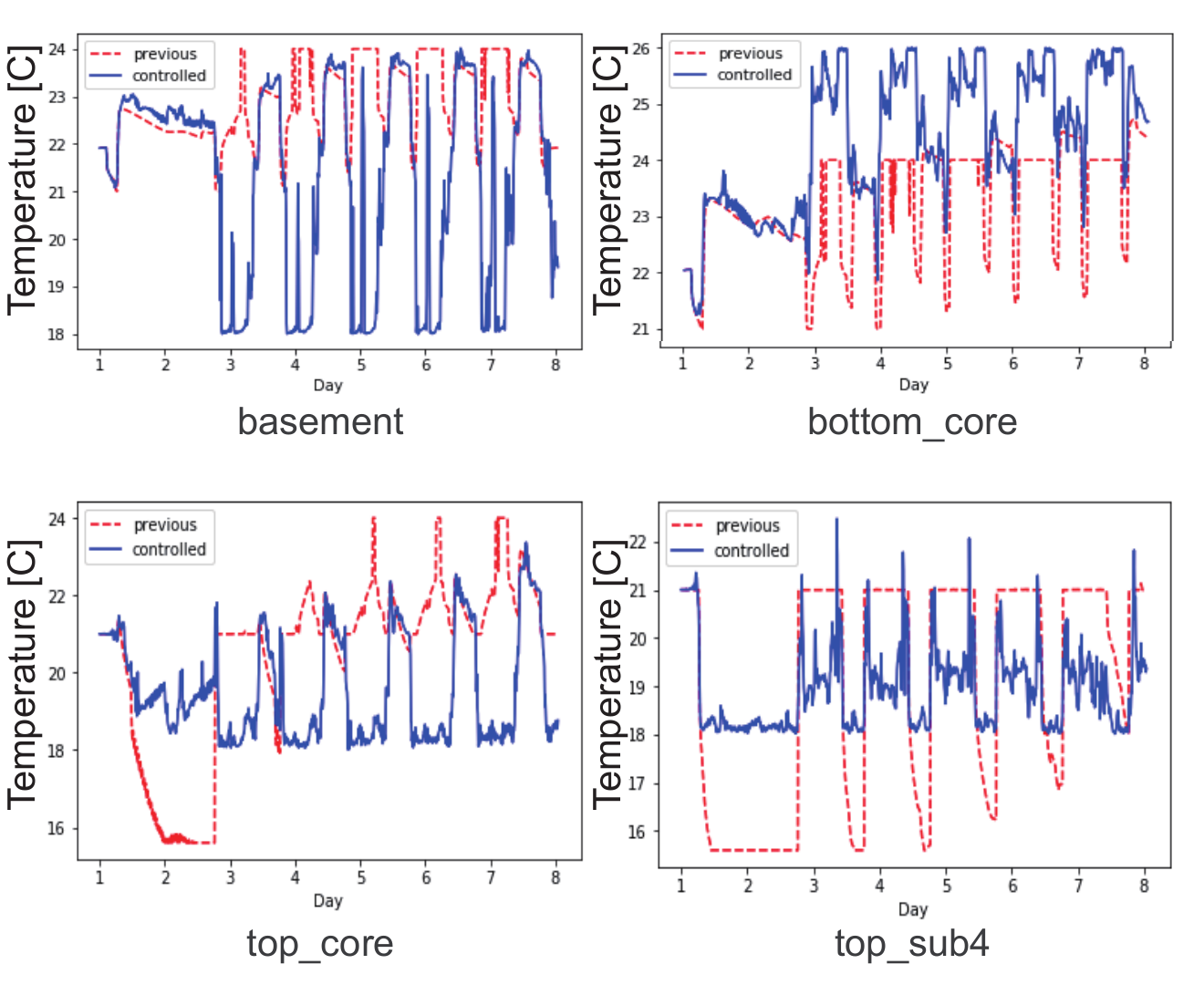}
	\caption{\small
		{
			One week temperature setpoint profile for 4 different zones of the office building.
		}
	}
	\label{zones}
\end{figure}

Fig.~\ref{zones} demonstrates how our proposed approach is able to find a group of control inputs for the building system globally. All of four zones' setpoint schedule exhibit daily patterns. Yet they are set to different values and evolution patterns. These setpoint schedule can provide to building operators, and it remains to be examined in real buildings if such optimized schedules could benefit the complex system as a whole.

